\documentclass[11pt]{amsart}
\usepackage{amssymb,amsmath}

\theoremstyle{plain}
\newtheorem{thm}{Theorem}[section]
\newtheorem{theorem}[thm]{Theorem}

\newtheorem{cor}[thm]{Corollary}

\theoremstyle{definition}

\newtheorem{ex}[thm]{Example}

\theoremstyle{remark}
\newtheorem{rem}[thm]{Remark}

\newtheorem*{remark*}{Remark}

\numberwithin{equation}{section}

\newcommand{\xb}{{\bf x}}

        \newcommand{\field}[1]{{\mathbb{#1}}}
        \newcommand{\NN}{\field{N}}

        \newcommand{\RR}{\field{R}}

\newcommand{\Tr}{\mbox{\rm Tr}}

\allowdisplaybreaks

\begin{document}
\title[Semiclassical spectral asymptotics
for a Schr\"odinger operator]{Semiclassical spectral asymptotics\\
for a magnetic Schr\"odinger operator\\ with non-vanishing magnetic
field}

\author{Bernard Helffer}

\address{D\'epartement de Math\'ematiques, B\^atiment 425, Univ
Paris-Sud 11, F-91405 Orsay C\'edex, France}
\email{Bernard.Helffer@math.u-psud.fr}

\author{Yuri A. Kordyukov }
\address{Institute of Mathematics, Russian Academy of Sciences, 112 Chernyshevsky
str. 450008 Ufa, Russia} \email{yurikor@matem.anrb.ru}

\thanks{B.H. is partially supported by INSMI
CNRS and by the ANR programme Nosevol. Y.K. is partially supported
by the Russian Foundation of Basic Research, projects 12-01-00519-a
and 13-01-91052-NCNI-a, and by the Ministry of education and science
of Russia, project 14.B37.21.0358. }

\begin{abstract}
We consider a magnetic Schr\"odinger operator $H^h$, depending on
the semiclassical parameter $h>0$, on a compact Riemannian manifold.
We assume that there is no electric field. We suppose that the
minimal value $b_0$ of the intensity of the magnetic field $b$ is
strictly positive. We give a survey of the results on asymptotic
behavior of the eigenvalues of the operator $H^h$ in the
semiclassical limit.
\end{abstract}

\date{\today}
 \maketitle


\section{Introduction}
Let $ M$ be a compact oriented manifold of dimension $n\geq 2$
(possibly with boundary). Let $g$ be a Riemannian metric and $\bf B$
a real-valued closed 2-form on $M$. Assume that $\bf B$ is exact and
choose a real-valued 1-form $\bf A$ on $M$ such that $d{\bf A} = \bf
B$. Thus, one has a natural mapping
\[
u\mapsto ih\,du+{\bf A}u
\]
from $C^\infty_c(M)$ to the space $\Omega^1_c(M)$ of smooth,
compactly supported one-forms on $M$. The Riemannian metric allows
to define scalar products in these spaces and consider the adjoint
operator
\[
(ih\,d+{\bf A})^* : \Omega^1_c(M)\to C^\infty_c(M)\,.
\]
A Schr\"odinger operator with magnetic potential $\bf A$ is defined
by the formula
\begin{equation}\label{defH}
H^h = (ih\,d+{\bf A})^* (ih\,d+{\bf A})\,.
\end{equation}
Here $h>0$ is a semiclassical parameter. If $M$ has non-empty
boundary, we will assume that the operator $H^h$ satisfies the
Dirichlet boundary conditions.

 From the geometric point of view, the 1-form $\bf A$ defines a
Hermitian connection $\nabla_{\bf A} = d-i{\bf A}$ on the trivial
complex line bundle $\mathcal L$ over $M$. The curvature of this
connection is $-i\bf B$. Then the operator $H^h$ is related with the
associated covariant (or Bochner) Laplacian
\[
H_{\bf A}=\nabla_{\bf A}^*\nabla_{\bf A}
\]
by the formula
\[
H^h=h^2(\,d-ih^{-1} {\bf A})^* (\,d-ih^{-1}{\bf A})= h^2 H_{h^{-1}
\bf A}\,.
\]
This formula shows, in particular, that the semiclassical limit
$h\to 0$ is clearly equivalent to the large magnetic field limit.

We choose local coordinates $\xb =(x_1,\ldots,x_n)$ on $M$. We write the
1-form $\bf A$ in the local coordinates as
\[
{\bf A}= \sum_{j=1}^nA_j(\xb )\,dx_j\,,
\]
the matrix of the Riemannian metric $g$ as
\[
g(\xb )=(g_{j\ell}(\xb))_{1\leq j,\ell\leq n}\,,
\]
and its inverse as
\[
g(\xb)^{-1}=(g^{j\ell}(\xb))_{1\leq j,\ell\leq n}\,.
\]
We denote the determinant of $g$ by:
 $$|g(\xb )|=\det(g(\xb ))\,.
 $$
 Then the magnetic field $\bf B$ is given
by the following formula
\[
{\bf B}=\sum_{j<k}B_{jk}\,dx_j\wedge dx_k\,, \quad
B_{jk}=\frac{\partial A_k}{\partial x_j}-\frac{\partial
A_j}{\partial x_k}\,.
\]
Moreover, the operator $H^h$ has in these coordinates the form
\[
H^h=\frac{1}{\sqrt{|g(\xb)|}}\sum_{1\leq j,\ell\leq n}\left(ih
\frac{\partial}{\partial x_j}+A_j(\xb)\right) \left[\sqrt{|g(\xb)|}
g^{j\ell}(\xb) \left(ih \frac{\partial}{\partial
x_\ell}+A_\ell(\xb )\right)\right].
\]
In the case when $M=\mathbb R^n$ is the flat Euclidean space, the
operator $H^h$ takes the form
\begin{equation}\label{e:Hh-flat}
H^h=\sum_{1\leq j\leq n}\left(hD_{x_j}- A_j(\xb)\right)^2,
\end{equation}
where, as usual, $D_{x_j} =\frac1i\frac{\partial }{\partial x_j}$,
$j=1,\ldots,n$.

When $n=2\,$, the magnetic two-form $\bf B$ is a volume form on $M$
and therefore can be identified with the function $b\in C^\infty(M)$
given by
\[
{\bf B}=b\,d\xb_g\,,
\]
where $d\xb_g$ denotes the Riemannian volume form $M$ associated with
$g$.

When $n=3\,$, the magnetic two-form $\bf B$ can be identified with a
magnetic vector field $\vec{b}$ by the Hodge star-operator. If $M$
is the Euclidean space $\RR^3$, we have
\begin{equation}\label{e:b-vec}
\vec{b}=(b_1,b_2,b_3) = \operatorname{curl} {\bf
A}=(B_{23},-B_{13},B_{12})\,,
\end{equation}
with the usual definition of $\operatorname{curl}$.

We are interested in asymptotic behavior of the spectrum of the
operator $H^h$ in the semiclassical limit. This problem was studied
in
\cite{FoHel1,miniwells,HM,HM01,HelSj7,Mat,MatUeki,Mont}
(see \cite{Fournais-Helffer:book,luminy, Ray4} for surveys including the case of problems with boundary).

After the pioneering works by Kato \cite{Kat}  and his school, the
starting reference for the spectral analysis of self-adjoint
realizations of the magnetic Schr\"odinger operator is the paper by
Avron-Herbst-Simon \cite{AHS} where the role of the module of the
magnetic field in the three-dimensional case appears for the first
time. Further investigations were inspired by R. Montgomery
\cite{Mont}, who was asking ``Can we hear the locus of the magnetic
field'' (by analogy with the celebrated question by M. Kac). In
\cite{Mont}, this question was studied for the two-dimensional
magnetic Schr\"odinger operator. Motivated by the question of R.
Montgomery, the first author and Mohamed in \cite{HM} investigated
the asymptotic behavior of the low-lying eigenvalues of the
Dirichlet realization of the magnetic Schr\"odinger operator in the
case when the magnetic field vanishes. This study was continued more
recently in \cite{PaKw,Qmath10,miniwells,Dombrowski-Raymond} (see
also \cite{luminy}). The case when the magnetic field never vanishes
was analyzed in detail for the Dirichlet realization in the
two-dimensional case in \cite{HM01} and more recently in
\cite{2Dcase,2DKarasev, RV}. Moreover, there is a big literature
devoted to the spectral analysis of the Neumann realization because
of its connection with problems in superconductivity (see
\cite{Fournais-Helffer:book} and the references therein). Finally,
we do not also give a complete description of the semi-classical
results obtained in the case when an electric potential $V$ is
creating the main localization and refer to \cite{HelSj7} and
\cite{DiSj} for a presentation and references therein.

The purpose of this paper is to give a survey of the results
obtained in the case when the magnetic field never vanishes. First,
we suppose that $M$ is two-dimensional.
Let
\begin{equation}
b_0=\min_{\xb\in M}|b(\xb )|\,.
\end{equation}
Remark that if $M$ is without boundary then we necessarily have
$b_0=0$, since
\[
\int_M b(\xb )d\xb_g=\int_Md{\bf A}=0\,.
\]
If we suppose that $M$ has a non-empty boundary and the operator
$H^h$ satisfies the Dirichlet boundary conditions, it was observed
by many authors \cite{Mont,Mall,Ue1,Ue2} (as the immediate
consequence of the Weitzenb\"ock-Bochner type identity and the
positivity of the square of a suitable Dirac operator) that, if $U$
is a domain in $M$, then, for any $u\in C^\infty_c(U)$, the
following estimate holds:
\begin{equation}\label{e:Mont}
\|(ih\,d+{\bf A})u\|^2_U\geq h  \int_U b |u|^2d\xb _g\,.
\end{equation}
In particular, for any $h>0\,$,
\begin{equation}
\lambda_0(H^h)\geq hb_0\,.
\end{equation}
In the case $M=\mathbb R^2$, this estimate follows from the formula
\[
h b(x) = -i [hD_{x_1} -A_1, hD_{x_2}-A_2]\,,
\]
which implies (after an integration by parts) that
\[
h \int b(\xb ) |u(\xb )|^2\,d\xb \leq  \|(hD_{x_1} - A_1) u\|^2 + \|(h D_{x_2} -
A_2) u\|^2\,.
\]

Due to this estimate, the function $h b$ can be considered in many
spectral problems   as an effective electric potential, that is, as
a magnetic analog of the electric potential $V$ in a Schr\"odinger
operator $-h^2\Delta+V$.

Any connected component of the minimum set
\begin{equation}\label{defU}
U = \{\xb \in M\,:\, b(\xb )= b_0\}
\end{equation}
can be understood as a magnetic well (attached to the given energy
$hb_0$). In particular, an asymptotic description of the spectrum
near the bottom strongly depends on the geometry of the magnetic wells and the
behavior of $b$ near them.

In higher dimensions, the role of magnetic potential is played by
the function  $\xb \mapsto h\cdot {\Tr}^+ (B(\xb ))$, which can be defined
in the following way. For any $\xb \in M$, denote by $B(\xb )$ the
anti-symmetric linear operator on the tangent space $T_\xb { M}$
associated with the 2-form $\bf B$:
\[
g_\xb(B(\xb)u,v)={\bf B}_\xb (u,v),\quad u,v\in T_\xb { M}.
\]
Recall that the intensity of the magnetic field is defined as
\[
{\Tr}^+ (B(\xb ))=\sum_{\substack{\lambda_j(\xb )>0\\ i\lambda_j(\xb )\in
\sigma(B(\xb)) }}\lambda_j(\xb )=\frac{1}{2}\Tr([B^*(\xb )\cdot
B(\xb )]^{1/2}).
\]
In the $(3D)$ case the only positive eigenvalue is $|B(x)|$ and we get
$$
{\Tr}^+ (B(\xb )) =|B(\xb )|\,.
$$
In the general case, we do not have the equivalent of \eqref{e:Mont} but only the weaker estimate \cite{HM}:
\begin{equation}\label{lb}
(h \inf_\xb {\Tr}^+ (B(\xb ))  - C h^\frac 54 )\int | u(\xb
)|^2\,d\xb   \leq \langle H^h u\,,\, u\rangle\,,\, \forall u \in
C^\infty_c (M)\,.
\end{equation}

When $U$ is not connected, the spectrum is essentially obtained by
analyzing (the union of) the spectra of Dirichlet Laplacians
attached to each component. This is true modulo exponentially small
errors. This corresponds to the so called magnetic tunneling.  We
will not focus on this question (which is widely open) and will more
emphasize on the presentation of the known semi-classical results in
dimension $2$ and $3$ which are purely magnetic first at the bottom
(Sections \ref{s2} and \ref{s3} in dimension $2$ and Section
\ref{s4} in dimension 3), secondly in Section \ref{s5} for excited
states in dimension $2$ where we present the newest contributions
(Helffer-Kordyukov and Raymond-Vu Ngoc)  but will give a few
examples in the last  section.

\section{Discrete wells in dimension $2$}\label{s2}
In this section, we will discuss the case of discrete wells. We
assume that:
\begin{equation}
b_0>0\,,
\end{equation}
and that  there exist a unique point $x_0$, which belongs to the interior of $M$,
  $k\in \NN$ and $C>0$ such that for all $x$ in some neighborhood
  of $x_0$ the estimates hold:
\begin{equation}
C^{-1}\, d(x,x_0)^2\leq b(x)-b_0 \leq C \, d(x,x_0)^2\,.
\end{equation}
We introduce:
\[
a={\rm Tr}\left(\frac12 {\rm Hess}\,b(x_0)\right)^{1/2}, \quad
d=\det \left(\frac12 {\rm Hess}\,b(x_0)\right)^{1/2}\,,
\]
and denote by $\lambda_0(H^h)\leq \lambda_1(H^h)\leq \lambda_2(H^h)\leq
\ldots$ the eigenvalues of the operator $H^h$ in $L^2(M)\,$.

\begin{thm}\label{t:mainbis}
Under current assumptions, for any $j\in \mathbb N$, there exists a
sequence $(\alpha_{j,\ell})_{\ell\in \mathbb N}$ with
\[
\alpha_{j,0} = b_0,\quad \alpha_{j,1}=0, \quad \alpha_{j,2}=
\frac{2d^{1/2}}{b_0}j+ \frac{a^2}{2b_0}\,,
\]
such that
\begin{equation}\label{e:mainbis}
\lambda_j(H^h) \sim h \sum_{\ell=0}^N \alpha_{j,\ell} h^{\frac \ell
2}\,.
\end{equation}
In other words, for any $N$, there exist $C_{j,N}>0$ and $h_{j,N}>0$
such that, for any $h\in (0,h_{j,N}]$,
\[
|\lambda_j(H^h) - h \sum_{\ell=0}^N \alpha_{j,\ell} h^{\frac \ell
2}|\leq C_{j,N} h^{\frac {N+3}{2}}\,.
\]
\end{thm}

In particular, we have for the groundstate
energy $\lambda_0(H^h)$ a two term asymptotics:
\[
\lambda_0(H^h) = hb_0 +h^2\frac{a^2}{2b_0}+\mathcal O (h^{5/2}),
\quad h\to 0,
\]
and the asymptotics of the splitting between the groundstate energy
and the first excited state~:
\[
\lambda_1(H^h)- \lambda_0(H^h) \sim h^2 \frac{2d^{\frac 12}}{b_0}\,.
\]

This theorem is proved in \cite{2Dcase}. A two-terms asymptotics for
the ground state energy in the flat case was previously obtained in
\cite{HM01}. Recent improvements by Helffer-Kordyukov
\cite{2DKarasev} and Raymond-Vu Ngoc \cite{RV} (see also
Section~\ref{s:revisited}) show that no odd powers of $h^{\frac12}$
actually occur in the flat case.  We believe that this
fact also holds in the general case of Riemannian manifold.

The proof of the upper bound is based on a construction of
approximate eigenfunctions for the operator $H^h$. More precisely,
we prove  in \cite{2Dcase}  the following accurate upper bound for the
eigenvalues of the operator $H^h$.

\begin{thm}
Under current assumptions, for any  $j$ and $k$ in $\mathbb N$, there exists
a sequence $(\mu_{j,k,\ell})_{\ell\in \mathbb N}$ with
\[
\mu_{j,k,0} = (2k+1)b_0, \quad \mu_{j,k,1}=0\,,
\]
and
\[ \mu_{j,k,2}= (2j+1)(2k+1)\frac{d^{1/2}}{b_0}
+(2k^2+2k+1)\frac{t}{2b_0} +\frac{1}{2}(k^2+k) R(x_0)\,,
\]
where $R$ is the scalar curvature, and
\[
t={\rm Tr}\left(\frac12 {\rm Hess}\,b(x_0)\right)\,,
\]
and for any $N$, there exist $\phi^h_{jkN}\in C^\infty(M)$,
$C_{jk,N}>0$ and $h_{jk,N}>0$ such that
\begin{equation}\label{e:orth}
(\phi^h_{j_1k_1N},\phi^h_{j_2k_2N}) =\delta_{j_1j_2}\delta_{k_1k_2}+
\mathcal O_{j_1,j_2,k_1,k_2}(h)\,,
\end{equation}
and, for any $h\in (0,h_{jk,N}]$,
\[
\|H^h\phi^h_{jkN}- \mu_{jkN}^h \phi^h_{jkN}\|\leq
C_{jkN}h^{\frac{N+3}{2}}\|\phi^h_{jkN}\|,\,
\]
where
\begin{equation}\label{e:mu}
\mu_{jkN}^h=h\sum_{\ell=0}^N \mu_{j,k,\ell} h^{\frac \ell 2}\,.
\end{equation}
\end{thm}

Since the operator $H^h$ is self-adjoint, using the Spectral Theorem, we
immediately deduce the existence of eigenvalues near the values
$\mu_{jkN}^h\,$.

\begin{cor}\label{c:dist}
For any  $j$, $k$ and $N$ in $\mathbb N$, there exist $C_{jk,N}>0$ and
$h_{jk,N}>0$ such that, for any $h\in (0,h_{jk,N})\,$,
\[
{\rm dist}(\mu_{jkN}^h, {\rm Spec}(H^h))\leq
C_{jk,N}h^{\frac{N+3}{2}}\,.
\]
\end{cor}

\begin{rem}
The low-lying eigenvalues of the operator $H^h\,$, as $h\rightarrow
0$, are obtained by taking $k=0$ in Theorem~\ref{t:qmodes}.
Therefore, as an immediate consequence of Theorem~\ref{t:qmodes}, we
deduce that, for any  $j$ and $N$ in $\mathbb N$, there exists $h_{j,N}>0$
such that, for any $h\in (0,h_{j,N}]\,$, we have
\[
\lambda_j(H^h)\leq \mu_{j0N}^h + C_{j0,N}h^{\frac{N+3}{2}}\,.
\]
In particular, this implies the upper bound in Theorem~\ref{t:mainbis}.
\end{rem}

\begin{rem}
Our interest in the case of arbitrary $k$ in Theorem~\ref{t:qmodes}
is motivated, in particular, by its importance for proving the
existence of gaps in the spectrum of the operator $H^h$ in the
semiclassical limit \cite{diff2006}.
\end{rem}

\begin{rem}\label{r:Landau}
The term
\[
(2k+1) hb_0 + \frac12   h^2\left(k^2+k\right)R\,,
\]
in the right-hand side of \eqref{e:mu} (see also \eqref{e:lambdakx}
below) has a natural interpretation as Landau levels. The
interpretation depends on whether $R$ is zero, positive or negative
and, in all three cases, is given in terms of eigenvalues of the
associated magnetic Laplacian with constant magnetic field (Landau
operator) on the corresponding simply connected Riemann surface of
constant curvature (see \cite{2Dcase} for more details).

We also mention the paper \cite{Ferapontov-Veselov} by Ferapontov and Veselov, who prove
that these three model magnetic Laplacians are integrable in some
sense. This observation enables them to give the complete description of
the spectra of these operators in the same way as it was done by
Schr\"odinger for the harmonic oscillator.
\end{rem}

\section{Degenerate wells in dimension $2$}\label{s3}
In this section, following
\cite{2Dcase2}, we will discuss the case when the minimum of the
magnetic field is attained on a regular curve $\gamma$. We assume that:
\begin{itemize}
  \item $b_0>0\,$;
  \item the set $\{x\in M : |b(x)| =b_0\}$ is a smooth curve $\gamma$,
  which is contained in the interior of $M$;
\item there is a constant $C>0$ such that for all $x$
in some neighborhood of $\gamma$ the estimates hold:
\begin{equation}\label{YK:B1}
C^{-1}d(x,\gamma)^2\leq |b(x)|-b_0 \leq C d(x,\gamma)^2\,.
\end{equation}
\end{itemize}
\subsection{Asymptotics  near the bottom}
The main purpose is to give  an asymptotics of the groundstate
energy $\lambda_0(H^h)$ of the operator $H^h$. Denote by $N$ the
external unit normal vector to $\gamma$. Let $\tilde{N}$ denote the
natural extension of $N$ to a smooth normalized vector field on $M$,
whose integral curves starting from a point $x$ in  a tubular
neighborhood of $\gamma$ are the minimal geodesics  to $\gamma$.
Consider the function $\beta_2$ on $\gamma$ given by
\begin{equation}\label{e:def-beta}
\beta_2(x)=\tilde{N}^2 |b(x)|\,, \quad x\in \gamma\,.
\end{equation}
By (\ref{YK:B1}), it is easy to see that
\[
\beta_2(x)>0\,, \quad x\in \gamma\,.
\]

\begin{theorem} \label{c:equiv-lambda0}
There exists $h_0>0$, such that, for any $h\in (0,h_0]\,$,
\begin{equation}
\lambda_0(H^h) =  hb_0 + h^2\,\frac{\mu_0}{4b_0} + \mathcal O (
h^{17/8}) \,.
\end{equation}
where
\begin{equation}\label{defmu0}
\mu_0:=\inf_{x\in\gamma} \beta_2(x)\,.
\end{equation}
\end{theorem}

The proof of the upper bound is based on a construction of
approximate eigenfunctions for the operator $H^h$. We denote by $R$ the
scalar curvature of the Riemannian manifold $(M,g)$.

\begin{theorem}\label{t:qmodes}
For any $x\in \gamma $ and for any integer $k\geq 0$, there exist
$C$ and $h_0>0$, such that, for any $h\in (0,h_0]$, there exists
$\Phi^h_k\in C^\infty_c(M), \Phi^h_k\neq 0\,,$ such that
\[
\left\|H^h\Phi^h_k-  \lambda^h(k,x)\Phi^h_k\right\|\leq C h^{17/8}
\|\Phi^h_k\|\,,
\]
where
\begin{equation}\label{e:lambdakx}
\lambda^h(k,x)=(2k+1) hb_0 +
h^2\left[(2k^2+2k+1)\frac{\beta_2(x)}{4b_0}+\frac12
\left(k^2+k\right) R(x)\right]\,.
\end{equation}
\end{theorem}

When $k=0$, we get:
\begin{cor}\label{Coroll}
For any $x\in \gamma $, there exist $C$ and $h_0>0$, such that, for
any $h\in (0,h_0]$, there exists $\Phi^h_0\in C^\infty_c(M),
\Phi^h_0\neq 0\,,$ such that
\[
\left\|H^h\Phi^h_0-\lambda^h(x)\Phi^h_0\right\|\leq C h^{17/8}
\|\Phi^h_0\|\,,
\]
where
\[
\lambda^h(x)= hb_0 + h^2\,\frac{\beta_2(x)}{4b_0}\,.
\]
\end{cor}

\subsection{Miniwells}\label{s:miniwells}
Like in the case of the Schr\"odinger operator with electric
potential (see \cite{HelSj5}), one can introduce an internal notion
of magnetic well for a fixed closed curve $\gamma$ in the minimum
set of the magnetic field $\mathbf B$. Such magnetic wells can be
naturally called magnetic miniwells. They are defined by means of
the function $\beta_2$ on $\gamma$ given by \eqref{e:def-beta}.

\begin{thm}
Assume that there exists a unique minimum point $x_0\in \gamma$ of
the function $\beta_2$ on $\gamma$, which is nondegenerate~:
\[
\mu_2:= \beta^{\prime\prime}_2(x_0)>0\,.
\]
For any  $j\in \mathbb N$, there exist $C_j$ and $h_j>0$, such that for
any $h\in (0,h_j)$
\[
\lambda_j(H^h)\leq hb_0 + h^2\,\frac{\mu_0}{4b_0}+
h^{5/2}\,\frac{(\mu_0\mu_2)^{1/2}}{4b_0^{3/2}}(2j+1) + C_j
h^{11/4}\,.
\]
\end{thm}

Here and below the derivative means the derivative with respect to
the natural parameter on $\gamma$.

\begin{rem}
We conjecture that
$$
\lambda_0(H^h)=  hb_0 + h^2\,\frac{\mu_0}{4b_0}+
h^{5/2}\,\frac{(\mu_0\mu_2)^{1/2}}{4b_0^{3/2}}+ o ( h^{5/2})\,.
$$
\end{rem}

The proof is based on a construction of approximate eigenfunctions,
which can be made near an arbitrary Landau level. For $k\in\NN$,
consider the function $V_k$ on $\gamma$ given by (cf.
\eqref{e:lambdakx})
\begin{equation}\label{e:defVk}
V_k(x):=(2k^2+2k+1)\frac{\beta_2(x)}{4b_0}+\frac12
\left(k^2+k\right) R(x)\,.
\end{equation}
Assume that there exists a unique minimum $x_0\in \gamma$ of the
function $V_k$ on $\gamma$, which is nondegenerate, that is
satisfying, for all $x\in \gamma$ in some neighborhood of $x_0\,$,
\begin{equation}\label{e:wells1}
Cd(x,x_0)^2\leq  V_k(x) - V_k(x_0)\leq C^{-1}d(x,x_0)^2\,.
\end{equation}
Under these assumptions, one can give the following, more precise
construction of approximate eigenvalues of the operator $H^h$.

\begin{thm}\label{t:qmodes0}
Under current assumptions, for any $j, k \in \mathbb N\,$, there exist
$u^h_{jk}\in C^\infty_c(M)$, $C_{jk}>0$ and $h_{jk}>0$ such that
\[
(u^h_{j_1k},u^h_{j_2k}) =\delta_{j_1j_2}+ \mathcal O_{j_1,j_2,k}(h)
\]
and, for any $h\in (0,h_{jk}]\,$,
\[
\|H^hu^h_{jk}- \mu_{jk}^h u^h_{jk}\|\leq
C_{jk}\, h^{11/4}\|u^h_{jk}\|\,,
\]
where
\begin{equation}\label{e:mujk}
\mu_{jk}^h= \mu_{j,k,0}h+\mu_{j,k,4} h^2+ \mu_{j,k,6} \, h^{5/2}\,,
\end{equation}
with
\[
\mu_{j,k,0} = (2k+1)b_0,\quad \mu_{j,k,4} = V_k(x_0)\,,
\]
and
\[
\mu_{j,k,6}=\frac{1}{2b_0}
V^{\prime\prime}_k(x_0)^{1/2}\beta_2(x_0)^{1/2}(2k+1)^{1/2}(2j+1)\,.
\]
\end{thm}

\section{Excited states for discrete wells}\label{s:revisited}\label{s4}
If Theorem \ref{t:mainbis} is satisfactory for the analysis of a
finite numbers of eigenvalues at the bottom, it appears to be useful
to get a extended description of the bottom of the spectrum
including more excited states. Motivated by Karasev's paper
\cite{Kar}, it seems to be interesting to produce an effective
Hamiltonian whose spectrum will also describes the excited states.
In 2013, Helffer-Kordyukov \cite{2DKarasev} on one side, and
Raymond-Vu Ngoc \cite{RV} on the other side reanalyze the problem in
the case of discrete wells with two different points of view leading
in the two cases to the existence of an effective $(1D)$-Hamiltonian
whose spectrum describes the spectrum of our magnetic Schr\"odinger
operator.

\subsection{Using a Grushin's problem (after Helffer-Kordyukov \cite{2DKarasev})}\label{ss41}
The approach of  \cite{2DKarasev} is based on Grushin's method. This
method was initiated in the context of hypoellipticity by V. Grushin
\cite{Gru} and then exploited by J.~Sj\"ostrand alone or with
collaborators in many contexts. We refer to \cite{SjZw} for a survey
on this method and references or Appendix D in \cite{Zw}. In
spectral theory a variant of this method is known under the name of
``Feschbach projection method'' or ``Schur complement formula'' in
analytic Fredholm theory.

Let us consider the magnetic Schr\"odinger operator $H^h$ in the
flat Euclidean space ${\mathbb R}^2\,$:
\[
H^h=h^2D_x^2+(hD_y+A(x,y))^2\,.
\]
The magnetic field $\bf B$ is given by
\[
{\bf B}=b\,dx\wedge dy\,\mbox{ with }\quad b(x,y) =\frac{\partial
A}{\partial x}(x,y)\,.
\]

Let
\[
b_0=\min_{(x,y)\in {\mathbb R}^2} |b(x,y)| >0 \,.
\]
We assume that at $\infty$, we have
$$
b_0 < \liminf_{|x|+|y| \rightarrow  +\infty} |b(x,y)|:=b_0+\eta_0
\,.
$$

Then one can prove easily \cite{HM} that, for any $0\leq
\eta_1<\eta_0 \,$, there exists $h_1>0$ such that
\[
\sigma (H^h)\cap [0, h(b_0+\eta_1)) \subset \sigma_d (H^h)\,, \quad
\forall h\in (0,h_1].
\]

Next, as above, we assume that:
\begin{itemize}
\item $b_0>0$;
\item the set $\{(x,y)\in {\mathbb R}^2 : |b(x,y)| =b_0\}$ is a single point $(x_0,y_0)$;
\item $(x_0,y_0)$ is a non-degenerate minimum:
\[
\operatorname{Hess} b (x_0,y_0)>0.
\]
\end{itemize}
We have a diffeomorphism $\phi : {\mathbb R}^2\to {\mathbb R}^2$
defined by
\[
\phi(x,y)=(A(x,y),y), \quad (x,y)\in {\mathbb R}^2\,.
\]
We then associate with $b$  a function $\hat b\in C^\infty({\mathbb
R}^2)$ by
\[
\hat b = b\circ \phi^{-1}\,.
\]

\begin{thm}\label{mainth}
There exist $h_0>0, \epsilon_0>0, \gamma_0 \in (0,\eta_0)$,
$h\mapsto \gamma_0(h)$ defined for
   $(0,h_0]$ such that $\gamma_0(h)\to \gamma_0$ as $h\to 0$,
and a semiclassical symbol $p_{\rm eff}(y, \eta ,h, z)$, which is
defined in a neighborhood $\Omega\subset {\mathbb R}^2$ of the set
$\{(y,\eta)\in \mathbb R^2: \hat b(y,\eta) \leq b_0+\gamma_0\}$ for
$h\in (0,h_0]$ and  $z\in \mathbb C$ such that
$|z|<\gamma_0+\epsilon_0$, of the form
\begin{equation}\label{bth}
p_{\rm eff} (y, \eta ,h, z) \sim\sum_{j\in \mathbb N} p_{\rm
eff}^{j}(y,\eta,z) h^j\,,
\end{equation}
with
\begin{equation}\label{cth}
 p_{\rm eff}^0 (y,\eta,z) =  \hat b(y,\eta) - b_0 -z\,,
\end{equation}
such that $\lambda_h \in \sigma (H^h)\cap [0, h (b_0
+\gamma_0(h)))$,  if and only if the associated
$h$-pseudo\-diffe\-ren\-tial operator\footnote{We use the Weyl
semi-classical quantization of the symbol (see for example
\cite{Hor})} $p_{\rm eff}(y, hD_y,h, z(h))$ has an approximate
$0$-eigen\-fun\-cti\-on $u_h^{qm}\in C^\infty(\mathbb R)$, i.e.
\begin{equation}\label{ath}
p_{\rm eff}(y, hD_y,h, z(h)) u_h^{qm} = \mathcal O (h^\infty)\,,
\end{equation}
with
\[
z(h)= \frac 1 h (\lambda_h - h b_0)+\mathcal O (h^\infty)\,,
\]
$|z(h)|<\gamma_0(h)$ for any $h\in (0,h_0]$, and such that the
frequency set\footnote{See \cite{Zw}  for a discussion of the frequency set and references therein.} of $u_h^{qm}$ is non-empty and contained in
$\Omega\,$.
\end{thm}
\begin{rem}
Here \eqref{ath} makes sense modulo $\mathcal O (h^\infty)$ by
extending first the symbol $p_{\rm eff}(y, \eta ,h, z)$ outside the
neighborhood $\Omega$ to a semiclassical symbol in $\mathbb R^2$ and
defining then the operators $p_{\rm eff}(y, hD_y,h, z)$ by the Weyl
calculus. Using the localization of the frequency set of $u_h^{qm}$,
the left hand side of \eqref{ath} does not depend on the extension
up to an error which is $\mathcal O (h^\infty)$.
\end{rem}
\begin{rem}
By Theorem~\ref{mainth}, for any $E\in [b_0,b_0+\gamma_0)$, the
spectrum of the operator $H^h$ (divided by $h$) is determined near
$E$ (say in an interval $(E-C h^\frac 12, E+ C h^\frac 12)$) and
modulo $\mathcal O(h^\frac 32)$ by the spectrum of $\hat b(y,hD_y) +
h b_1(y, hD_y,E)$, where one can use the Bohr-Sommerfeld rule (see
\cite{HelRo} or \cite{HS88}  for a mathematical justification) for
determining the energy levels.
\end{rem}

\begin{cor} There exists $\gamma_0 \in (0, \eta_0)$, $h_0>0$ and $C>0$ such that
$$
\lambda_{j+1}(H^h) -\lambda_j(H^h) \geq \frac 1C h^2 \,,\, \forall
h\in (0,h_0]\,,
$$
for any $j$ such that $\lambda_{j+1}(H^h) < h (b_0 + \gamma_0)$.
\end{cor}

\subsection{Using a Birkhoff Normal form (after Raymond--Vu Ngoc
\cite{RV})}\label{s:RV} The proof of Raymond--Vu Ngoc is reminiscent
of Ivrii's approach (see his book (old version or new version in
progress on his Home Page) \cite{Iv1} and the -- more accessible but
without proofs -- introductory article \cite{Iv0}) and uses a
Birkhoff normal form. This approach has the advantage to be
semi-global and uses more general symplectomorphisms and their
quantizations.

Consider the magnetic Schr\"odinger operator $H^h$ in $\mathbb R^2$
given by \eqref{e:Hh-flat}. Let $H$ be its $h$-symbol:
\begin{multline}\label{defsymbH}
H(x,y,\xi,\eta)=|\xi-A_1(x,y)|^2+|\eta-A_2(x,y)|^2, \\
(x,y,\xi,\eta)\in T^*\mathbb R^2=\mathbb R^2\times \mathbb R^2\,.
\end{multline}
By definition the energy surface $\Sigma_E$ corresponding to energy $E$ is the set $H^{-1} (E)$.
The first result shows the existence of a smooth symplectic
diffeomorphism that transforms the initial Hamiltonian into a normal
form, up to any order in the distance to the zero energy surface $\Sigma_0$.
Assume that the magnetic field $b$ does not vanish in an open set
$\Omega\subset \mathbb R^2$.

\begin{thm}[\cite{RV}, Theorem 1.1] \label{th4.5}
There exists a symplectic diffeomorphism $\Phi$, defined in an open
set $\tilde \Omega \subset \mathbb C_{z_1}\times \mathbb R^2_{z_2}$,
with values in $T^*\mathbb R^2$, which sends the plane $\{z_1=0\}$
to $\Sigma_0$, and such that
\[
H\circ \Phi=|z_1|^2f(z_2,|z_1|^2)+\mathcal O(|z_1|^\infty),
\]
where $f:\mathbb R^2\times \mathbb R\to \mathbb R$ is smooth.
Moreover, the map
\[
\varphi : \Omega \ni (x,y)\mapsto \Phi^{-1}(x,y,\mathbf A(x,y))\in
(\{0\}\times \mathbb R^2_{z_2})\cap \tilde \Omega
\]
is a local diffeomorphism and
\[
f\circ (\varphi(x,y),0)=|b(x,y)|.
\]
\end{thm}

The next result gives the quantum counterpart of this theorem. We keep the notation of the previous theorem.

\begin{thm}[\cite{RV}, Theorem 1.6]
For $h$ small enough there exists a (semi-classical) Fourier
Integral Operator\footnote{See \cite{Hor,Zw} for a definition.}
$U_h$ such that
\[
U^*_hU_h=I+Z_h, \quad U_hU_h^*=I+Z^\prime_h,
\]
where $Z_h, Z^\prime_h$ are $h$-pseudo-differential operators that
microlocally vanish in a neighborhood of $\tilde \Omega\cap \Sigma_0$,
and
\[
U^*_hH^hU_h=\mathcal I_hF_h+R_h,
\]
where:
\begin{enumerate}
  \item $\mathcal I_h:=-h^2\frac{\partial^2}{\partial x_1^2}+x_1^2$\,.
  \item $F_h$ is a classical $h$-pseudo-differential operator
  that commutes with $\mathcal I_h$\,.
  \item For any Hermite function $h_n(x_1)$ such that $\mathcal
  I_hh_n=h(2n-1)h_n$, the operator $F^{(n)}_h$ acting on
  $L^2(\mathbb R_{x_2})$ by
  \[
h_n\otimes F^{(n)}_h(u)=F_h(h_n\otimes u)
  \]
  is a classical $h$-pseudo-differential operator with principal symbol
  \[
F^{(n)}(x_2,\xi_2)=b(x,y)\,,
  \]
  where $(0,x_2+i\xi_2)=\varphi(x,y)\,$.
  \item Given any $h$-pseudo-differential operator $D_h$ with principal
  symbol $d_0$ such that $d_0(z_1,z_2)=c(z_2)|z_1|^2+\mathcal
  O(|z_1|^3)$, and any $N\geq 1$, there exist classical pseudo-differential
  operators $S_{h.N}$ and $K_N$ such that
  \[
R_h=S_{h.N}(D_h)^N+K_N+\mathcal O(h^\infty)\,,
  \]
  with $K_N$ compactly supported away from a fixed neighborhood of
  $|z_1|=0$.
  \item $\mathcal I_hF_h=\mathcal N_h=\mathcal H^0_h+Q_h\,$, where
  $\mathcal H^0_h$ is the $h$-pseudodifferential operator of symbol
  $H^0(z_1,z_2)=b(\varphi^{-1}(z_2))|z_1|^2\,$, and the operator $Q_h$
  is relatively bounded with respect to $\mathcal H^0_h$ with an
  arbitrarily small relative bound.
\end{enumerate}
\end{thm}

As a consequence, Raymond and Vu Ngoc obtain the following theorem.

\begin{thm}[\cite{RV}, Theorem 1.5] \label{t:RV-Theorem 1.5}
Assume that the magnetic field $B$ is non vanishing on $\mathbb R^2$
and confining: there exist constants $\tilde C_1>0$, $M_0>0$ such
that
\[
b(q)\geq \tilde C_1 \ for \ |q|\geq M_0\,.
\]
Let $\mathcal H^0_h=Op^w_h(H^0)$, where
$H^0=b(\varphi^{-1}(z_2)|z_2|^2$ where $\varphi :\mathbb R^2\to
\mathbb R^2$ is a diffeomorphism. Then there exists a bounded
classical pseudo-differential operator $Q_h$ on $\mathbb R^2$, such
that
\begin{itemize}
  \item $Q_h$ commutes with $Op^w_h(|z_1|^2)$;
  \item $Q_h$ is relatively bounded with respect to $\mathcal H^0_h$
  with an arbitrarily small relative bound;
  \item its Weyl symbol is $O_{z_2}(h^2+h|z_1|^2+|z_1|^4)\,$,
\end{itemize}
so that the following holds. Let $0<C_1<\tilde C_1$. Then the
spectra of $H^h$ and $\mathcal N_h:=\mathcal H^0_h+Q_h$ in
$(-\infty, C_1h]$ are discrete. We denote by $0<\lambda_1(h)\leq
\lambda_2(h)\leq\cdots $ the eigenvalues of $H^h$ and by
$0<\mu_1(h)\leq \mu_2(h)\leq\cdots $ the eigenvalues of $\mathcal
N_h$. Then for any $j\in \mathbb N^*$ such that $\lambda_j(h)\leq
C_1h$ and $\mu_j(h)\leq C_1h$, we have
\[
|\lambda_j(h)-\mu_j(h)|=O(h^\infty)\,.
\]
\end{thm}

\begin{rem}
Theorem \ref{t:RV-Theorem 1.5} is stronger than in
Theorem~\ref{mainth} because Theorem~\ref{mainth} gives a
description of the spectrum of $H^h$ in the interval
$[hb_0,h(b_0+\gamma_0))$ for some $\gamma_0\in (0,\eta_0)$, whereas
in Theorem \ref{t:RV-Theorem 1.5}, $\gamma_0\in (0,\eta_0)$ is
arbitrary. On the other hand, the symbol of the effective
Hamiltonian in Theorem~\ref{t:RV-Theorem 1.5} seems to be less
explicit than in Theorem~\ref{mainth}. The other point could be that
Theorem~\ref{mainth} allows us to treat an additional term $h^2
V(x,y)$. This will complete the analysis of Helffer-Sj\"ostrand
\cite{HSLNP345}, in the case of the constant magnetic field. The
case with an additional term $h V$ could also be interesting.
\end{rem}

\begin{rem}
As communicated to us by F. Faure, there is some hope that the
results of \cite{RV} can be generalized under a generic assumption
to the case of arbitrary even dimension. Some results are also
presented in \cite[Chapter 13]{Iv1}.
\end{rem}

\section{Discrete wells in dimension 3}\label{s5}
In this section, we discuss the three-dimensional case.
\subsection{Upper bounds
\cite{3Dcase}} Consider the magnetic Schr\"o\-din\-ger operator
$H^h$ in a domain $\Omega$ of the flat Euclidean space ${\mathbb
R}^3$ (see \eqref{e:Hh-flat}).  As usual, we assume that $H^h$
satisfies the Dirichlet boundary condition. Let
$\vec{b}=(b_1,b_2,b_3)$ be the corresponding vector magnetic field
(see \eqref{e:b-vec}).

We assume that there exists a constant $C>0$ such that for $j=1,2,3$
we have
\begin{equation}\label{ass0}
|(\nabla b_j)(\xb )|\leq C(|\vec{b}(\xb )|+1), \quad \forall \xb \in
\Omega\,.
\end{equation}
Put
\[
b_0=\min \{|\vec b(\xb )|\, :\, \xb \in \Omega\}.
\]
We assume that there exist a (connected) bounded domain
$\Omega_1\subset\subset \Omega$ and a constant $\epsilon_0>0$ such
that
\begin{equation}\label{ass1}
 |\vec b(\xb )| \geq b_0+\epsilon_0, \quad \xb \not\in \Omega_1\,.
\end{equation}
As shown in \cite{HM}, under conditions \eqref{ass0} and
\eqref{ass1}, for any $\epsilon_1$ with $0<\epsilon_1<\epsilon_0$,
there exists $h_1>0$ such that, for $h\in (0,h_1]$
\[
\sigma(H^h)\cap [0,h(b_0+\epsilon_1))\subset \sigma_d(H^h).
\]
Denote by $\lambda_0(H^h)\leq \lambda_1(H^h)\leq \lambda_2(H^h)\leq
\ldots$ the eigenvalues of the operator $H^h$ contained in
$[0,h(b_0+\epsilon_0))$.

Finally, we assume that:
\[
b_0>0\,,
\]
and that there exists a unique minimum $\xb_0\in \Omega$ such that
$|\vec b(\xb_0)|=b_0$, which is non-degenerate: in some neighborhood
of $\xb_0$
\[
C^{-1}|\xb-\xb_0|^2\leq |\vec b(\xb)|-b_0 \leq C |\xb-\xb_0|^2\,.
\]
We also introduce:
\[
d=\det {\rm Hess}\,|\vec b|(\xb_0)\,, \quad a=\frac{1}{2 b_0^2}({\rm
Hess} |\vec b|\, {\vec b}\cdot \vec b) (\xb_0)\,.
\]

\begin{thm}\label{t:main3}
Under current assumptions, for any  $m\in \mathbb N$, there exist $C_m>0$
and $h_m>0$ such that, for any $h\in (0,h_m]\,$,
\begin{equation}\label{exp3D}
\lambda_m(H^h) \leq h b_0 +h^{3/2}
a^{1/2}+h^2\left[\frac{1}{2b_0}\left(\frac{d}{2a}\right)^{1/2}(2m+1)+
\nu \right]+C_m h^{9/4}\,,
\end{equation}
where $\nu $ is some explicit constant\footnote{which means that it is given
by a rather complicated explicit formula}.
\end{thm}

The proof of Theorem \ref{t:main3} is based on a construction of
quasimodes.
\begin{thm}\label{t:qmodes3}
Under current assumptions, for any  $j$, $k$ and $m$ in $\mathbb N$, there
exist $\phi^h_{j,k,m}\in C^\infty_c(\Omega)$, $C_{j,k,m}>0$ and
$h_{j,k,m}>0$ such that
\[
(\phi^h_{j_1,k_1,m_1},\phi^h_{j_2,k_2,m_2})
=\delta_{j_1j_2}\delta_{k_1k_2}\delta_{m_1m_2}+ \mathcal
O_{j_1,k_1,m_1,j_2,k_2,m_2}(h)\,,
\]
and, for any $h\in (0,h_{j,k,m}]$,
\[
\|H^h\phi^h_{j,k,m}- \mu_{j,k,m}^h \phi^h_{j,k,m}\|\leq C_{j,k,m}\,
h^{\frac{9}{4}}\|\phi^h_{j,k,m}\|,\,
\]
where
\[
\mu_{j,k,m}^h=\mu_{j,k,m,0}h+\mu_{j,k,m,2} h^{\frac 3
2}+\mu_{j,k,m,4} h^2\,.
\]
with
\[
\mu_{j,k,m,0} = (2k+1)b_0\,,\quad
\mu_{j,k,m,2}=(2j+1)(2k+1)^{1/2}a^{1/2}\,,
\]
and
\[
\mu_{j,k,m,4}=\frac{1}{2b_0}\left(\frac{d}{2a}\right)^{1/2}(2m+1)(2k+1)+
\nu(j,k)\,,
\]
where $\nu(j,k)$ has the form
\[
\nu(j,k)=\nu_{22}(2k+1)^2+\nu_{11}(2j+1)^2+\nu_0\,,
\]
with some explicit constants $\nu_0\,,\,\nu_{11},\,\nu_{22}\,$.
\end{thm}

\begin{rem}
It is conjectured that
\[
\lambda_m(H^h) \geq hb_0
+h^{3/2}a^{1/2}+h^2\left[\frac{1}{2b_0}\left(\frac{d}{2a}\right)^{1/2}(2m+1)+
\nu \right]- C_m h^{9/4}\,.
\]
At the moment, we only know from  \cite{HM}
\[
\lambda_m(H^h) \geq hb_0 - C h^2\,,
\]
which is an improvement of the general lower bound \eqref{lb}.
\end{rem}

\subsection{On some statements of V. Ivrii \cite{Iv0,Iv1}}
Here we refer to some results announced in \cite{Iv0} and developed
in Chapter 18 in \cite{Iv1}. These results correspond in the
$(3D)$-case  to what was discussed in the $(2D)$-case in the
subsection~\ref{s:RV}. Under the assumption that the magnetic field
does not vanish, the claim\footnote{We have tried to correct many
typos of the statement (more specifically the remainder in Formula
(25) ) in \cite{Iv0}. Note in particular that the sum in the right
hand side is undefined (see however \cite{2DKarasev} which meets the
same problem) and could only be meaningful for some subspace of functions
whose energy is for example less than $h b_0 + C h^\frac 32$.} is
that (up to conjugation by an $h$-Fourier integral operator), our
Schr\"odinger operator can microlocally be written in the form:
\begin{multline*}
\omega_1 (x_1,x_2, h D_{x_2}) (h^2 D_{x_3}^2  + x_3^2) + h^2 D_{x_1}^2 \\
 + \sum_{2m +n  +\ell\geq 3} h^\ell
\,a_{mn\ell} (x_1,x_2, h D_{x_2})\, (h^2 D_{x_3}^2+ x_3^2)^m (h
D_{x_1})^n\,,
\end{multline*}
 with
 $$
 \omega_1 = |\vec b| \circ \psi\,,
 $$
where $\psi$ is some local unspecified diffeomorphism which plays
the role of $\varphi^{-1}$ in Theorem \ref{th4.5}.

Once precisely stated and proved, let us explain what we could
expect after. Reducing to the lowest Landau level (the first
eigenvalue of $h^2 D_{x_3}^2 + x_3^2$), we obtain that the spectrum
of our initial operator near the minimum of $|\vec b(x)|$ should be
deduced from the spectral analysis in $(-\infty, hb_0 + C h^\frac
32)$ (for some fixed $C>0$) in the semi-classical limit of the
following ``formal'' pseudo-differential operator:
\[
h \omega_1 (x_1,x_2, h D_{x_2}) +  h^2 D_{x_1}^2  + \sum_{2m+n +\ell
\geq 3} h^{m+\ell}\, a_{mn\ell} (x_1,x_2, h D_{x_2}) (h
D_{x_1})^n\,.
\]
If we only look for the principal term (and divide by $h$), we get
as first ``effective'' operator to analyze for the spectrum now in
$(-\infty, b_0+ Ch^\frac 12)$:
\[
\omega_1(x_1,x_2,h D_{x_2}) + (h^\frac 12 D_{x_1})^2 + h a_{020} (x_1,x_2, hD_{x_2})\,,
\]
with the hope to get in this way an approximation modulo $\mathcal O
(h^\frac 54)$. This suggests a semi-classical analysis near the
bottom of a pseudodifferential operator of the type met in
Born-Oppenheimer theory $p(x_1,x_2, h^\frac 12 D_{x_1}, h D_{x_2})$
with two semi-classical parameters (see \cite{Ray4} and references
therein for a recent discussion on this subject). This would be
coherent with the expansion obtained in the right hand side of
\eqref{exp3D} at least modulo $\mathcal O (h^\frac 94)\,$.

\section{Some remarks and open questions}
\subsection{Geometry of magnetic fields} Consider the
magnetic Schr\"o\-din\-ger operator in the flat Euclidean space
${\mathbb R}^n$ (see \eqref{e:Hh-flat}). Its semiclassical symbol
(as defined in \eqref{defsymbH}) is a smooth function $H\in
C^\infty(\mathbb R^{2n})$ whose zero set of $H$ given by
\[
\Sigma_0:=H^{-1}(0)=\{(\xb,\xi) \in \mathbb R^{2n} : \xi_j=A_j(\xb)\,,
j=1,\ldots,n \}\,.
\]
Since it is a graph, it is an embedded submanifold of $\mathbb
R^{2n}$, parameterized by $\xb \in \mathbb R^{n}$. It is easy to check
that if we denote by $J: \mathbb R^{n}\to \Sigma$ the embedding
$J(\xb)=(\xb ,A(\xb ))$, then, for the canonical symplectic form $\omega
=\sum_{j=1}^nd\xi_j\wedge dx_j$ on $\mathbb R^{2n}$ we have
\[
J^*\omega\left|_\Sigma\right.\cong \mathbf B\,.
\]

When $n=2$ and the magnetic field $b$ does not vanish, $\Sigma_0$ is
symplectic. When $n=3$, $\Sigma_0$ cannot be symplectic. If the
magnetic field $\vec b$ does not vanish, then $\mathbf B$ has
constant rank, and $\Sigma_0$ is a presymplectic manifold. When $n$
is arbitrary even, we can hope that, under generic assumptions,
$\Sigma_0$ is symplectic. This kind of analysis was basic in the
seventies for the analysis of the hypoellipticity of operators with
multiple characteristics.

Recall that, in Remark~\ref{r:Landau}, we give a geometric
interpretation of some terms, entering into the asymptotic formula
\eqref{e:mu} for approximate eigenvalues of the operator $H^h$ in
the two-dimensional case. One can naturally consider similar
questions in the three-dimensional case. First, observe that three
cases $R=0$, $R>0$ and $R<0$ mentioned in Remark~\ref{r:Landau}
correspond to three cases of two-dimensional model geometries:
Euclidean, spherical and hyperbolical, respectively. In the
three-dimensional case, the situation is more complicated. There are
eight three-dimensional model geometries introduced by Thurston
(see, for instance, \cite{Thurston}). The interesting open problem
is to construct the magnetic Schr\"odinger operators with constant
magnetic field on each three-dimensional geometric model and compute
its spectra. It is also interesting to find examples of integrable
magnetic Schr\"odinger operators on three-dimensional Riemannian
manifolds.

\subsection{The tunneling effect} Although, as a
consequence of magnetic Agmon estimates \cite{HelSj7, HM, Ray4}, it
is possible to give upper bounds on the tunneling effect (see
Section 7.2 in \cite{Zw} or \cite{He1} for an introduction) due to
the presence of multiconnected magnetic wells, essentially no
results are  known for lower bounds of this effect analogous to what
is proved for the celebrated double well problem for the
Schr\"odinger operator $-h^2 \Delta + V$. The only exception is
\cite{HelSj7}, which involves $\sum_j (hD_{x_j} - t(h) A_j)^2 + V$
but this last result is not a ``pure magnetic effect'' and it is
assumed that the magnetic field is small enough ($|t(h)| = \mathcal
O (h |\log h|)$).

There are however a few models where one can ``observe'' this effect
in particular in domains with corners \cite{BDMV} (numerics with
some theoretical interpretation, see also
\cite{Fournais-Helffer:book} for a presentation of results due to V.
Bonnaillie-Noel), the role of the magnetic wells being played by the
corners of smallest angle. We describe other toy models, which are
closer to the analysis which is presented in this survey:

\begin{ex}
We consider in $\mathbb R^2$ the operator:
$$ h^2 D_{x}^2 + (h D_y-a(x))^2 + y^2\,.$$ This
model is rather artificial (and not purely magnetic)  but by Fourier
transform, it is unitary equivalent to
$$
h^2 D_x^2 + (\eta -a(x))^2 + h^2 D_\eta^2\,,
$$
which can be analyzed because it enters in the category of the
miniwells problem treated in Helffer-Sj\"ostrand  \cite{HelSj5}. We
have indeed a well defined in $\mathbb R^2_{x,\eta}$  by $\eta =
a(x)$ which is unbounded but if we assume a varying curvature $
\beta (x)= a'(x)$ (with $ \liminf_{|x| \rightarrow+\infty} |\beta
(x)| > \inf_x |\beta (x)|$) we will have a miniwell localization. A
double well phenomenon can be created by assuming $\beta =a'$ even.
\end{ex}

\begin{ex}
If we add an electric potential $V(x)$ to the previous example, we
get:
$$ h^2 D_{x}^2 + (h D_y-a(x))^2 + y^2 + V(x) \,.$$
For $a(x)=x$, this example was considered by J. Br\"uning, S. Yu.
Dobrokhotov and R.V.\,Nekrasov in \cite{BDN}.

Here one can measure the explicit effect of the magnetic field by considering
$$
h^2 D_{x}^2 + h^2 D_\eta^2 + (\eta -a(x))^2 + V(x)\,.
$$
If $V$ admits as minimum value $0$, the wells are defined in
$\mathbb R^2_{x,\eta}$ by $\eta = a(x), V(x)=0$ and one can use
under suitable assumptions the semi-classical treatment of the
double well problem for the Schr\"odinger operator with electric
potential $W(x,\eta)=(\eta -a(x))^2 + V(x)$ (see \cite{He1}).
\end{ex}

\begin{ex}
One can also imagine that in the case of Sections~\ref{s2} and \ref{s4},  we have a
magnetic double well, and that a tunneling effect could be measured
using the effective $(1D)$-hamiltonian introduced in Subsection \ref{ss41} $\hat b (x,hD_x)$ (actually a perturbation of it), assuming
that $b$ and $A$ are holomorphic with respect to one of the variables. Here
we are extremely far for a proof but we could hope for candidates
for a formula for the splitting.
\end{ex}

\begin{ex}
Similarly, one can hope to measure the tunneling in the case of
miniwells, in the situation considered in Subsection
\ref{s:miniwells}, when $|b|$ admits its minimum along a curve and
$\beta_2$ has two symmetric miniwells.
\end{ex}

\begin{ex}
 Finally one can come back to the Montgomery example
\cite{Mont} which was analyzed in \cite{HM,PaKw,Qmath10,
Dombrowski-Raymond} and  corresponds to the two dimensional case
when the magnetic field vanish to some order on a compact curve.
According to a personal communication of V. Bonnaillie-No\"el, F.
H\'erau and N. Raymond, it seems to be reasonable to hope (work in
progress) that one could analyze  the splitting between the two lowest eigenvalues for the following model in $\mathbb
R^2$:
\[
h^2D_x^2 + \left(hD_y - \gamma(y) \frac{x^2}{2}\right)^2\,,
\]
where $\gamma$ is a positive even $C^\infty$ function with two
non degenerate minima and $\inf \gamma < \liminf \gamma$. By
dilation, this problem is unitary equivalent to the analysis of the
spectrum of
$$
h^\frac 43 \left( D_x^2 + \left(h^\frac 13 D_y - \gamma(y)
\frac{x^2}{2}\right)^2\right)\,.
$$
After division by $h^\frac 43$, the guess is then that we can
understand the tunneling by analyzing the spectrum of the $h^\frac
23$-pseudodifferential operator on $L^2(\mathbb R)$ whose Weyl
symbol is $\gamma(x)^{\frac 23} E (\gamma(x)^{-\frac 13} \xi)\,$,
where $E(\alpha)$ is the ground state energy of the Montgomery
operator $D_t^2 + (\frac{t^2}{2}-\alpha)^2$. This would involve a
Born-Oppenheimer analysis like in \cite{Ray4}.
\end{ex}

\subsection*{Acknowledgment}
We thank V. Bonnaillie-No\"el, S. Dobrokhotov, F. Faure, F. H\'erau,
M. Karasev,  N. Raymond, S. Vu Ngoc for useful discussions or
exchange of information.

\end{document}